\documentclass[leqno,12pt]{article}
\usepackage{latexsym}
\usepackage{amsmath}
\usepackage{amssymb}
\usepackage{bbm}
\usepackage{amsthm}
\usepackage{pst-all}

\usepackage{graphicx}
\usepackage{epsfig}

 0 \DeclareMathAlphabet{\mathbfit}{OT1}{cmr}{bx}{it}

\newcommand{\bea}{\begin{eqnarray*}}
\newcommand{\eea}{\end{eqnarray*}}
\newcommand{\be}{\begin{eqnarray}}
\newcommand{\ee}{\end{eqnarray}}
\newcommand{\beq}{\begin{equation}}
\newcommand{\eeq}{\end{equation}}

\newcommand{\zd}{\mathbb{Z}^d}

\newtheorem{thm}{Theorem}

\newtheorem{lem}[thm]{Lemma}

\newtheorem{quest}[thm]{Question}

\begin{document}

\title{Reconstructing the environment seen by a RWRE} 

\author{
{\small Nina Gantert}\\
\and
{\small Jan Nagel}\\
}

\maketitle

\begin{abstract}
  \noindent Consider a walker performing a random walk in an i.i.d. random environment, and assume that the walker tells us at each time the environment it sees at its present location. Given this history of the transition probabilities seen from the walker - but not its trajectory - can we tell if the RWRE is recurrent or transient? Can we reconstruct the law of the environment? We show that in a one-dimensional environment, the law of the environment can be reconstructed, and we know in particular if the RWRE is recurrent or transient. 
  
\end{abstract}

\medskip


\section{Introduction}

For a fixed mapping $\omega:\mathbb{Z}\rightarrow (0,1)$ the random walk $X:\mathbb{N}_0\rightarrow \mathbb{Z}$ in the environment $\omega$ is the Markov chain starting in 0 and law $P_\omega$ given by the transition probabilities 
\begin{align*}
P_\omega\big( X(n+1)=x+1|\, X(n)=x\big) &= \omega(x) ,\\
P_\omega\big( X(n+1)=x-1|\, X(n)=x\big) &= 1-\omega(x) .
\end{align*}
If we endow the set $\Omega$ of all environments with a probability measure $P$, this process is called a \emph{random walk in random environment (RWRE)} and $P_\omega$ is called \emph{quenched law}. We assume that $P=\mu^{\otimes \mathbb{Z}}$ is a product measure with marginal $\mu$.
We refer to \cite{So} and \cite{Z} for results on the RWRE, for instance, a criterion for recurrence and transience; our arguments will not need them.

In this paper we deal with the following question: Suppose that we only observe the sequence
\begin{align*}
\xi := \big( \xi(0),\xi(1),\dots \big) := \big( \omega(X(0)),\omega(X(1)),\dots \big) ,
\end{align*}
the history of transition probabilities at the walker's position, but we do not know the trajectory, is it possible to recover the marginal $\mu$? 

Of course, the same question may be asked for RWRE on $\zd$. Let ${\cal P}^d$ denote the set of probability measures on 
$\{+e_i, -e_i, 1 \leq i \leq d\}$ where $e_1, \ldots ,e_d$ are the unit vectors of $\zd$. For a fixed mapping
 $\omega:\zd \rightarrow {\cal P}^d$, the random walk $X:\mathbb{N}_0\rightarrow \zd$ in the environment $\omega$ is the Markov chain starting in 0 and law $P_\omega$ given by the transition probabilities 
\begin{align*}
P_\omega\big( X(n+1)=x\pm e_i|\, X(n)=x\big) &= \omega(x, \pm e_i) .
\end{align*}
Suppose that we only observe the sequence
\begin{equation}\label{observations}
\xi := \big( \xi(0),\xi(1),\dots \big) := \big( \omega(X(0), \cdot),\omega(X(1), \cdot),\dots \big) ,
\end{equation}
the history of transition probabilities at the walker's position, but we do not know the trajectory, is it possible to recover the marginal $\mu$? Note that in contrast to the one-dimensional case, recovering $\mu$ will not always tell us if the RWRE is recurrent or transient - despite some recent progress, there is still no criterion for recurrence/transience of RWRE in an i.i.d environment on $\zd$.

These questions are motivated from the classical scenery reconstruction problem, where the walk is a simple random walk on $\zd$, $ d \in \{1,2\}$,
the ``scenery'' is a colouring of $\mathbb{Z^d}$ and the walker tells us at each time the colour of its present location.
Given this sequence of observations - but not the trajectory of the walker - can we then reconstruct the scenery (up to translations, reflections and rotations)?
This problem goes back to Itai Benjamini and Harry Kesten, see \cite{BenKe}, and has lead to lot of interesting research, we 
refer to \cite{Kestensurvey} for some nice (and still open!) problems. One direction of research is to ask about the ergodic properties of the observation sequence given by (\ref{observations}), see \cite{FrankJeffsurvey}.
In the one-dimensional case, the original question can be answered in the positive in the sense that an i.i.d. scenery can be reconstructed almost surely, up to translation and reflection, see \cite{Henry2colours}.
In the two-dimensional case, this is possible if the number of colours in the scenery is large enough, see \cite{HenryMatthias}.
Typically, scenery reconstruction is easier if the scenery has more colours. Clearly, if $ d \in \{1,2\}$ and each location has a different colour, the trajectory of the walk can be reconstructed from the sequence of observations (up to symmetries, i.e. reflections and rotations).
In the same way, in our case, the question will be much easier if $\mu$ has a non-atomic part, cf. the argument below.

A related, but different question for RWRE was asked by Omer Adelman and Nathana\"el
Enriquez in \cite{AdelmanEn}: if we know a single ``typical'' trajectory of the walk, can we reconstruct the law of the environment? This questions is answered in the positive by \cite{AdelmanEn} for i.i.d. environments on $\zd$.

If $d=1$, the reconstruction of $\mu$ is possible, which is made more precise in the following theorem. We denote by $\mathcal{M}$ the set of all probability measures on $(0,1)$ and endow it with the weak topology and the corresponding Borel $\sigma$-algebra.

\medskip

\begin{thm}\label{haupt1}
Assume $d=1$. There exists a measurable mapping $\mathcal{A}:(0,1)^{\mathbb{N}_0}\rightarrow \mathcal{M}$, such that for any measure $\mu \in \mathcal{M}$
\begin{align*}
P_\omega \big( \mathcal{A}(\xi) = \mu \big) =1
\end{align*}
for $P$-almost all $\omega$.
\end{thm}

\medskip

Before we begin the proof, let us consider the simple case, where $\mu$ has a non-atomic part. We denote by $\mu=\mu_a+\mu_{na}$ the decomposition into the atomic part $\mu_a$ and the non-atomic part $\mu_{na}$. Note that whether $\mu_{na}$ is non-zero can be read off from the observations: the support of $\mu$ is almost surely equal to the closure of $\mathcal{S} = \{\xi(0),\xi(1),\dots \}$ and the set of atoms is almost surely given by  
\begin{align*}
\mathcal{S}_a = \{ \eta \in \mathcal{S} |\, \exists k\geq 0:\, \xi(k)=\xi(k+1)=\eta \} ,
\end{align*}
as only atoms can appear twice in the environment. Now $\mu_{na} $ is non-zero if and only if $\mathcal{S}_{na}=\mathcal{S}\setminus \mathcal{S}_a \neq \emptyset$. In this case, observations $\xi(n)\in \mathcal{S}_{na}$ can be used as perfect markers, as $\xi(n)=\xi(m)$ implies $X(n)=X(m)$. We give an informal description how this allows a reconstruction of $\mu$:

\begin{itemize}
\item Wait for observations $(\dots,\xi(m-2),\xi(m-1),\xi(m),\dots )$ in $\xi$, where $\xi(m-2),\xi(m-1) \in \mathcal{S}_{na}$, $\xi(m)\neq \xi(m-2)$ and both ``markers'' $\xi(m-2)$ and $\xi(m-1)$ have never before appeared in the sequence of observations.
\item If the assumptions are met for the $n$-th time, denote by $\eta_n$ the value of the corresponding $\xi(m)$.
\item When the two markers are seen for the first time, $X(m)$ must be at a point not visited before and $\xi(m)$ is the value of the environment at a point not visited before. Also, the choice of $\xi(m)$ is independent of the earlier entries of $\xi$, which implies $P_\omega$-almost surely
\begin{align*}
\tfrac{1}{n} \sum_{k=1}^n \delta_{\eta_k} \xrightarrow[n \rightarrow \infty ]{w} \mu .
\end{align*}
\end{itemize}

The perfect markers in $\xi$ immediately reveal whether $X$ is recurrent or transient. In the recurrent case, the following procedure constructs (a.s.) an environment which is up to translation equal to $\omega$.
\begin{itemize}
\item Choose two values $\eta_1,\eta_2\in \mathcal{S}_{na}$. 
\item Among all words $(\xi(n),\xi(n+1),\dots ,\xi(n+m))$ in $\xi$ with $\xi(n)=\eta_1, \xi(n+m)=\eta_2$ (of which there are infinitely many), there will be infinitely many of minimal length $m$. This word corresponds to a straight path of $X$ from $X(n)$ to $X(m)$. Therefore, $(\xi(n),\dots ,\xi(n+m))$ is a block of transition probabilities appearing in this order (or reversed) in $\omega$.
\item Repeat this step with new end points $\eta_2,\eta_3$ with a new marker $\eta_3 \in \mathcal{S}_{na}$ and concatenate the two obtained blocks of transition probabilities. It may happen that one block has to be contained in the other, which is the case if $\eta_1$ appears in the second block or $\eta_3$ in the first.
\item Continuing in this way, we obtain in the limit an environment $\hat\omega$, which is up to translation either $\omega$ or $\tilde\omega$, where $\tilde\omega(z) = \omega(-z)$ is the reflected environment. 
\item To decide on the orientation, consider all movements from a point $z$ with $\hat\omega_z\neq \tfrac{1}{2}$. A proportion of $\hat\omega_z$ of those movements needs to be made to the right.
\end{itemize}

\section{The main idea: the environment as a random walk} 
We now assume that $\mu$ is a purely atomic measure. We follow \cite{Henry3colours}. There can only be countably many support points $\eta_1,\eta_2,\dots$ which can be found in $\mathcal{S}$. We denote by $N \in [1,\infty]$ the cardinality of $\mathcal{S}$ and exclude the deterministic environments where $N=1$. Let $\mathcal{T}$ be the rooted tree with root $o$ where each vertex has exactly $N$ neighbours. We label the vertices by a mapping $\varphi:\mathcal{T}\rightarrow \mathcal{S}$ which satisfies
\begin{itemize}
\item $\varphi(o) = \xi(0)$
\item $\varphi$ restricted to the neighbours of any vertex $v$ is a bijection. That is, each vertex has exactly one neighbour which is labeled by a specific $\eta_i$.
\end{itemize}
Given an environment $\omega$, we define $R:\mathbb{Z}\rightarrow \mathcal{T}$ to be the bi-infinite path on $\mathcal{T}$ with $R(0)=o$ and $\varphi ( R(z)) = \omega(z)$ for all $z\in \mathbb{Z}$. Due to the second property in the definition of $\varphi$, this determines $R$ uniquely. As the environment is chosen under $P$ in an i.i.d.\ way, $R$ performs under $P$ a random walk on $\mathcal{T}$, starting at the root. In each step, both on the positive and negative time axis, $R$ moves from a vertex $v$ to a neighbour $w$ with probability $\mu(\varphi(w))$. Roughly speaking, this provides us with an embedding of the environment into the tree. Note that since we do not observe $\omega$, we do not know the path of $R$.\\ 
In a second step, the random walk $X$ on $\mathbb{Z}$ can be represented as a random walk $T$ on the trajectory of $R$. Given $X$, there is exactly one $T:\mathbb{Z} \rightarrow \{\dots, R(-1),R(0),R(1),\dots \}$ such that
\begin{align*}  
T(n) = (R\circ X)(n)
\end{align*}
for all $n\in \mathbb{N}_0$. The crucial point is that although we observe neither $X$ nor $R$, we know the path of $T$, as the labels of vertices visited by $T$ must coincide with the observation $\xi$, we have
\begin{align*}
(\varphi \circ T)(n) = (\varphi \circ R\circ X)(n) = (\omega \circ X)(n) = \xi(n) ,
\end{align*}
which we observe. In other words, as $X$ performs a random walk on $\mathbb{Z}$ and yields $\xi$, the process $T$ moves along a path on the tree giving the same observation $\xi$ when reading the labels of the vertices provided by $\varphi$. Due to the structure of the labeling, there is only one such path.\\

{\bf Example 1}\\
To illustrate this construction, we look at the case $N=2$, where the tree reduces to $\mathbb{Z}$ and the labeling by $\varphi$ is periodic repeating the word $\eta_0\eta_0\eta_1\eta_1$. Let us assume that the environment from position 0 to position 10 takes the values
\begin{align*}
(\omega(0),\dots ,\omega(10))=(\eta_0,\eta_0,\eta_1,\eta_0,\eta_1,\eta_1,\eta_0,\eta_0,\eta_0,\eta_1,\eta_0).
\end{align*}
The first observation at time 0 will be given by $\eta_0$, so we choose our labeling $\varphi$ such that $\varphi(0)=\varphi(1)=\eta_0$. This determines $\varphi$ uniquely on the whole integer line. The steps of $R$ representing this part of the environment are given by
\begin{align*}
(R(0),\dots ,R(10)) = (0,1,2,1,2,3,4,5,4,3,4) ,
\end{align*}
see Figure \ref{RundT} for an illustration. To keep this example simple, we do not consider $R$ in negative time, which corresponds to $\omega$ on the negative integers. Next, say the first steps of $X$ are as follows:
\begin{align*}
(X(0),\dots ,X(10)) = (0,1,2,3,4,3,4,5,6,7,6)
\end{align*}
This path gives us the observations
\begin{align*}
(\xi(0),\dots ,\xi(10)) = (\eta_0,\eta_0,\eta_1,\eta_0,\eta_1,\eta_0,\eta_1,\eta_1,\eta_0,\eta_0,\eta_0),
\end{align*}
which, given our choice of $\varphi$, implies the following movement of $T$:
\begin{align*}
(T(0),\dots ,T(10)) = (0,1,2,1,2,1,2,3,4,5,4)
\end{align*}
\\

\medskip

\psset{xunit=1cm,yunit=0.75cm,runit=0.5cm}
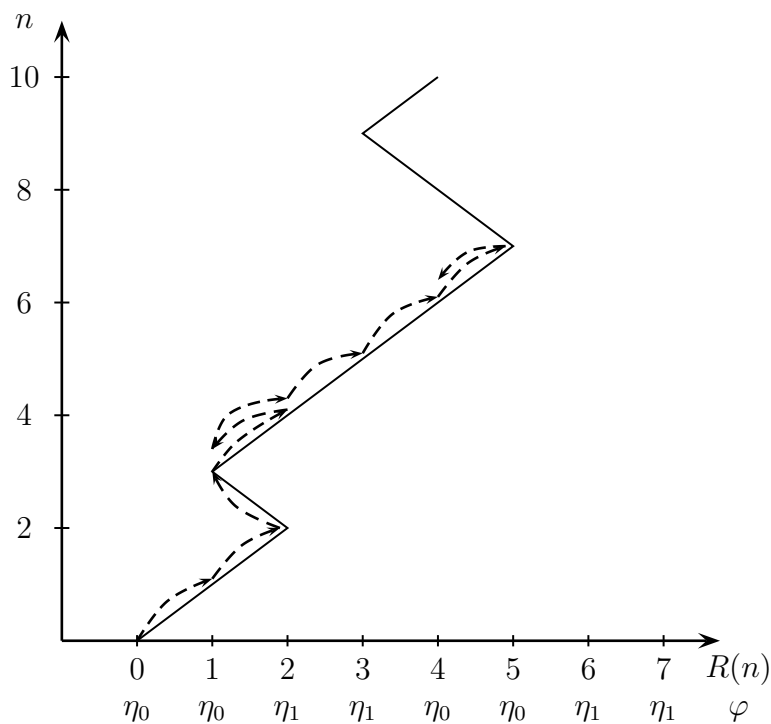
\begin{figure}[h] \label{RundT}
\centering
\begin{pspicture}(0,0)(12,14)

\psline[arrowsize=3pt 3,arrows=->,linewidth=1pt](1,2)(9.75,2)
\psline[arrowsize=3pt 3,arrows=->,linewidth=1pt](1,2)(1,13)

\psline[linewidth=0.75pt](2,1.9)(2,2.1)
\psline[linewidth=0.75pt](3,1.9)(3,2.1)
\psline[linewidth=0.75pt](4,1.9)(4,2.1)
\psline[linewidth=0.75pt](5,1.9)(5,2.1)
\psline[linewidth=0.75pt](6,1.9)(6,2.1)
\psline[linewidth=0.75pt](7,1.9)(7,2.1)
\psline[linewidth=0.75pt](8,1.9)(8,2.1)
\psline[linewidth=0.75pt](9,1.9)(9,2.1)

\psline[linewidth=0.75pt](0.9,4)(1.1,4)
\psline[linewidth=0.75pt](0.9,6)(1.1,6)
\psline[linewidth=0.75pt](0.9,8)(1.1,8)
\psline[linewidth=0.75pt](0.9,10)(1.1,10)
\psline[linewidth=0.75pt](0.9,12)(1.1,12)

\rput[cb](2,1.5){0}
\rput[cb](3,1.5){1}
\rput[cb](4,1.5){2}
\rput[cb](5,1.5){3}
\rput[cb](6,1.5){4}
\rput[cb](7,1.5){5}
\rput[cb](8,1.5){6}
\rput[cb](9,1.5){7}
\rput[cb](10,1.5){$R(n)$}
\rput[cb](2,0.75){$\eta_0$}
\rput[cb](3,0.75){$\eta_0$}
\rput[cb](4,0.75){$\eta_1$}
\rput[cb](5,0.75){$\eta_1$}
\rput[cb](6,0.75){$\eta_0$}
\rput[cb](7,0.75){$\eta_0$}
\rput[cb](8,0.75){$\eta_1$}
\rput[cb](9,0.75){$\eta_1$}
\rput[cb](10,0.75){$\varphi$}
\rput[cb](0.5,4){2}
\rput[cb](0.5,6){4}
\rput[cb](0.5,8){6}
\rput[cb](0.5,10){8}
\rput[cb](0.5,12){10}
\rput[cb](0.5,13){$n$}

\psline[linewidth=0.75pt](2,2)(4,4)(3,5)(7,9)(5,11)(6,12)

\pscurve[linewidth=1pt,linestyle=dashed]{->}(2,2)(2.4,2.7)(3,3.1)
\pscurve[linewidth=1pt,linestyle=dashed]{->}(3,3.1)(3.4,3.7)(3.9,4)
\pscurve[linewidth=1pt,linestyle=dashed]{->}(3.9,4)(3.3,4.4)(3,5)
\pscurve[linewidth=1pt,linestyle=dashed]{->}(3,5)(3.3,5.5)(4,6.1)
\pscurve[linewidth=1pt,linestyle=dashed]{->}(4,6.1)(3.4,5.9)(3,5.4)
\pscurve[linewidth=1pt,linestyle=dashed]{->}(3,5.4)(3.2,6)(4,6.3)
\pscurve[linewidth=1pt,linestyle=dashed]{->}(4,6.3)(4.4,6.9)(5,7.1)
\pscurve[linewidth=1pt,linestyle=dashed]{->}(5,7.1)(5.4,7.8)(6,8.1)
\pscurve[linewidth=1pt,linestyle=dashed]{->}(6,8.1)(6.4,8.7)(6.9,9)
\pscurve[linewidth=1pt,linestyle=dashed]{->}(6.9,9)(6.4,8.9)(6,8.4)

\end{pspicture}
\caption{The first moves of $R$ representing the enviroment and of $X$ as a random walk on the trajectory of $R$. The dashed arrows indicate the movements of $(R\circ X,X)$, the path of $T$ is obtained by projecting onto the first coordinate.}
\end{figure}

The process $T$ can only be transient if both $X$ and $R$ are transient, otherwise it is recurrent. Even though the increments of $R$ are not i.i.d.\ under $P$, the behaviour is essentially the same as for the simple random walk on the tree. 

\begin{lem} \label{recctran}
$R$ visits the root infinitely often if and only if $N=2$.
\end{lem}

In order to make statements about the movement of $X$ when we only observe $T$, we look for specific crossings of finite paths by $T$. For a generic process $S:I\rightarrow W$ with $I\subset \mathbb{Z}$ and $W$ a tree, we call $(i_1,i_2)$ a \emph{crossing} of $(w_1,w_2)$ by $S$, when $S(i_1)=w_1,S(i_2)= w_2$ and $S(i)\notin \{w_1,w_2\}$ for $\min\{i_1,i_2\}<i<\max \{i_1,i_2\}$. We call this crossing \emph{positive}, if $i_1<i_2$ and \emph{negative} otherwise. The crossing is said to be \emph{straight}, if $|i_2-i_1|$ is equal to the path distance between $w_1$ and $w_2$. \\
Consider again the example above, where $(0,5)$ is a crossing of $(0,3)$ by $R$. Since $R$ steps back during the time interval $(0,5)$, this is not a straight crossing. On the other hand, $(4,7)$ is a straight crossing of $(2,5)$ by $R$.\\ 
Of central importance to us are straight crossings of a path $(v_1,v_2)$ in the tree by $T$, as $T$ can only move in a straight way on the trajectory of $R$ if $R$ moves in a straight way on the tree $\mathcal{T}$. \\

\parbox[c]{0.9\textwidth}{If  $(i_1,i_2)$ is a straight crossing of $(v_1,v_2)$ by $T$ then $(i_1,i_2)$
is a straight crossing by $X$ of a straight crossing by $R$, that is, there are $(z_1,z_2)$ such that $(i_1,i_2)$ is a straight crossing of $(z_1,z_2)$ by $X$ and $(z_1,z_2)$ is a straight crossing of $(v_1,v_2)$ by $R$.}\hfill ($\ast$)\\

In our example, $(6,9)$ is a straight crossing of $(2,5)$ by $T$. Indeed, during the time $(6,9)$, $X$ performs a straight crossing of $(4,7)$ and $(4,7)$ is a straight crossing of $(2,5)$ by $R$, see figure 1.

\section{Proofs}

\medskip

\textbf{Proof of Theorem \ref{haupt1}:} \\
We first consider $N=2$, that is, $\mu =\lambda_0 \delta_{\eta_0} + \lambda_1 \delta_{\eta_1}$ with $\lambda_1=1-\lambda_0$ and $\mathcal{T}=\mathbb{Z}$. Without loss of generality, we assume $\xi(0)=\eta_0$ and choose our labeling $\varphi$ such that 
$\varphi(0)=\varphi(1)=\eta_0$. Consequently, we have $\varphi(4m)=\varphi(4m+1)=\eta_0$ and $\varphi(4m+2)=\varphi(4m+3)=\eta_1$ for all $m\in \mathbb{Z}$. For a stochastic process $Z$ on a tree, we denote by 
\begin{align*}
\tau_Z(v) = \inf\{n | \, Z(n)=v\}
\end{align*}
the hitting time of vertex $v$. For $m\geq 0$, define $I_m =(4m+1,4m+4)$ and let $W_m$ be the indicator random variable which is 1 if the first crossing of $I_m$ by $T$ is straight and 0 otherwise. By ($\ast$), $W_m=W_m^RW_m^X$, where $W_m^R$ is an indicator variable equal to 1 if and only if the first crossing of $I_m$ by $R$ is straight and $W_m^X$ is 1 if and only if the first crossing by $X$ of the first crossing of $I_m$ by R is straight. We will show that $W_0,W_1,\dots $ are independent and identically distributed, this time following \cite{Henry2colours}. \par

For the independence, note that conditioned on $R$ (or on $\omega$, i.e.\ under the quenched law $P_\omega$), the random variables $W_0^X, W_1^X,\dots $ depend only on the path of $X$ between ladder times of $X$ -- the times when $X$ reaches a point $z \in \mathbb{Z}$ where a crossing of a new $I_m$ by $R$ begins. That is, if $z_{1,m}=\tau_R(4m+1),z_{2,m}=\tau_R(4m+4)$, then $W_m^X$ depends only on
\begin{align*}
X(\tau_X(z_{1,m})+1)-X(\tau_X(z_{1,m})),X(\tau_X(z_{1,m})+2)-X(\tau_X(z_{1,m})),\dots ,X(\tau_X(z_{2,m}))-X(\tau_X(z_{1,m})) 
\end{align*}
and these collections of increments of $X$ are independent for different $m$, since $[z_{1,m},z_{2,m}]$ are disjoint intervals. This implies that $W_0^X, W_1^X,\dots $ are independent conditioned on $R$. Moreover, the conditional probability of the event $\{W_m^X=1\}$ depends only on the path segment of $R$ in the time interval $[z_{1,m},z_{2,m}]$ given by the random variable
\begin{align*}
R_m=\big( R(z_{1,m}+1)-R(z_{1,m}),\dots ,R(z_{2,m})-R(z_{1,m})\big),
\end{align*}
which again are independent for different $m$. In particular, $W_0^R,W_1^R,\dots $ are independent. Note that although $X$ may leave the corresponding path segment of $R$ during $[\tau_X(z_{1,m}),\tau_X(z_{2,m})]$, this does not influence the distribution of $W_m^R$. Consequently, we have
\begin{align*}
& \qquad P\big( W_{i_1}=1,\dots ,W_{i_k}=1 \big) \\
& =E\big[ P\big( W_{i_1}=1,\dots ,W_{i_k}=1 |\, R\big) \big] \\
& =E\big[ P\big( W^X_{i_1}=1,\dots ,W^X_{i_k}=1|\, R\big) \cdot \mathbbm{1}_{\{W^R_{i_1}=1,\dots ,W^R_{i_k}=1\}} \big]\\
& =E\big[ P\big( W^X_{i_1}=1|\, R\big)\cdots P\big( W^X_{i_k}=1|\, R\big) \cdot \mathbbm{1}_{\{W^R_{i_1}=1,\dots ,W^R_{i_k}=1\}} \big]\\
& =E\big[ P\big( W^X_{i_1}=1|\, R_{i_1}\big)\cdots P\big( W^X_{i_k}=1|\, R_{i_k}\big) \cdot \mathbbm{1}_{\{W^R_{i_1}=1,\dots ,W^R_{i_k}=1\}} \big]\\
& =E\big[ P\big( W^X_{i_1}=1|\, R_{i_1}\big) \mathbbm{1}_{\{W^R_{i_1}=1 \}} \big] \cdots 
E\big[ P\big( W^X_{i_k}=1|\, R_{i_k}\big) \mathbbm{1}_{\{W^R_{i_k}=1 \}} \big]\\
&= P\big( W^X_{i_1}=1 \big| W^R_{i_1}=1 \big)P\big( W^R_{i_1}=1 \big) \cdots P\big( W^X_{i_k}=1 \big| W^R_{i_k}=1 \big)  P\big( W^R_{i_k}=1 \big) \\
&= P\big( W_{i_1}=1 \big) \cdots P\big( W_{i_k}=1 \big) ,
\end{align*}
which proves the independence. \par 

We now evaluate the probability 
\begin{align*}
P(W_m=1) = P(W^X_m=1|\, W_m^R=1) P(W_m^R=1) .
\end{align*}
By our definition, $I_m$ is labeled as $(\eta_0,\eta_1,\eta_1,\eta_0)$ and $R$ moves to a neighbour labeled by $\eta_i$ with probability $\lambda_i$. Let $P_m$ denote the law of $R$ when starting at the left end point $4m+1$ and let $E_m$ be the event that $R$ reaches $4m+4$ before returning to $4m+1$. In order to reach $4m+4$ before returning to $4m+1$, $R$ needs to make two steps to the right (with probability $\lambda_1^2$), then make any number $k$ of steps between $4m+3$ and $4m+2$ and back before moving to $4m+4$. This gives 
\begin{align*}
P(W_m^R=1)= P_m(\tau_R(4m+4)=3|\, E_m) = \frac{P_m(\tau_R(4m+4)=3)}{P_m(E_m)} 
= \frac{\lambda_1^2\lambda_0}{\lambda_1^2 \left(\sum_{k=0}^\infty \lambda_1^{2k}\right) \lambda_0} = 1-\lambda_1^2 . 
\end{align*}
Given that the first crossing of $I_m$ by $R$ is straight, the probability of $\{W_m^X=1\}$ depends on whether $X$ moves on the positive or on the negative integers. If the first crossing of $I_m$ by $R$ happens during an interval 
$(t_{m,1},t_{m,2})$ in positive time (which corresponds to the environment on the positive integers), the process $X$ needs to move to the right for $T$ to cross $I_m$. In this case the first crossing of $I_m$ by $T$ is a crossing by $X$ of a positive crossing by $R$. If on the other hand the first crossing of $I_m$ by $R$ is by the trajectory in negative time, $X$ performs a crossing of the crossing by $R$ by moving to the left and the corresponding crossing of $R$ is negative. Let $D_m$ be the event that $t_{m,1}>0$, then
\begin{align*}
P(W_m^X=1|\, W_m^R=1, D_m) = \frac{\eta_0 \eta_1^2}{\eta_0\eta_1 \left(\sum_{k=0}^\infty ((1-\eta_1)\eta_1)^k\right) \eta_1}  = 1-(1-\eta_1)\eta_1,
\end{align*}
as $X$ moves from $t_{m,1}$ (or $t_{m,1}+1$) to $t_{m,1}+1$ (or $t_{m,1}+2$) with probability $\eta_0$ ($\eta_1$, respectively) and in the other direction with probability $1-\eta_0$ (and $1-\eta_1$). Given $D_m^c$, we need to interchange $\eta_i$ and $1-\eta_i$, which by our choice of $I_m$ leads to
\begin{align*}
P(W_m^X=1|\, W_m^R=1, D_m^c)  = 1-\eta_1(1-\eta_1) = P(W_m^X=1|\, W_m^R=1, D_m) .
\end{align*}
Using that $\{W_m^R=1\}$ is independent of $D_m$, we get $P(W^X_m=1|\, W_m^R=1) = 1-\eta_1(1-\eta_1)$ and therefore,
\begin{align*}
P(W_m=1) = \big( 1-\eta_1(1-\eta_1) \big)\big( 1-\lambda_1^2\big) .
\end{align*}
This proves that $W_0^X, W_1^X,\dots $ are independent identically Bernoulli-distributed random variables. By the law of large numbers, we have $P_\omega$-almost surely
\begin{align*}
\frac{1}{n(1-\eta_1(1-\eta_1))} \sum_{k=1}^n W_k \xrightarrow[n \rightarrow \infty ]{} 1-\lambda_1^2 .
\end{align*} 
Since the $W_m$ are functions of $\xi$, this convergence provides us with a measurable mapping which, given $\xi$ yields $\lambda_1$. In the case $N=2$, this already determines the measure $\mu$.\par

In the general case $N\geq 2$ we reduce this to the procedure above. Fix two values $\eta_0,\eta_1\in \mathcal{S}$ to which $\mu$ assigns weights $\lambda_0$ and $\lambda_1$. The intervals $I_m$ are now replaced by disjoint vertex-sets $(I_m(\eta_0,\eta_1))_{m\geq 0}$ in the tree $\mathcal{T}$, such that each set $I_m(\eta_0,\eta_1)$ contains exactly four neighbouring vertices $v_{1,m},\dots ,v_{4,m}$ with strictly increasing distance from the root and labels $\varphi(v_{1,m}) = \varphi(v_{4,m}) = \eta_0$ and $\varphi(v_{2,m})= \varphi(v_{3,m})=\eta_1$. When $T$ crosses the $m$-th of such a set for the first time without leaving this set of vertices, let $W_m(\eta_0,\eta_1)$ be equal to 1 if this crossing is straight and 0 otherwise. The same arguments as in the case $N=2$, this time conditioning on a movement of $R$ within $I_m(\eta_0,\eta_1)$, show that $W_0(\eta_0,\eta_1), W_1(\eta_0,\eta_1),\dots $ are again independent and 
\begin{align*}
P\big( W_m(\eta_0,\eta_1)=1\big) = \big( 1-\eta_1(1-\eta_1) \big)\big( 1-\lambda_1^2\big) .
\end{align*}
The law of large numbers allows us to recover $\lambda_1$ and repeating this with different choices of values $\eta_1$ shows that we can recover any $\lambda_i$. The (countable) combination of all these operations yields a weight vector $(\lambda_0,\lambda_1,\dots)$ as a measurable function of $\xi$ by which we can define $\mathcal{A}(\xi) = \sum_{k=0}^N \lambda_k \delta_{\eta_k}$.
\qed

\medskip

\textbf{Proof of Lemma \ref{recctran}:} \\
Suppose $N=2$, then the tree $\mathcal{T}$ is just $\mathbb{Z}$ and $\mu = \lambda_0\delta_{\omega_0} + \lambda_1\delta_{\omega_1}$. Without loss of generality, assume that $\varphi(0)=\varphi(1)=\eta_0$. We show that $R$, when only observed at points $4m,m\in \mathbb{Z}$, behaves as a symmetric random walk. Let $\tau_0=0$ and for $n\geq 0$,
\begin{align*}
\tau_{n+1} = \inf \{ k\geq \tau_n|\, X_k \in 4\mathbb{Z} \} .
\end{align*}
To move from $4m$ to $4m+4$ without backtracking to $4m$, $R$ needs to make two steps to the right, then any number $l$ of steps from $4m+2$ to $4m+1$ and back and any number $r$ of steps from $4m+2$ to $4m+3$ and back, and then move two steps further to the right. Similar to the calculations in the prove of Theorem \ref{haupt1}, we get
\begin{align*}
P(X(\tau_{n+1})=4m+4|\, X(\tau_{n})=4m) = \lambda_0\lambda_1 \left( \sum_{l,r\geq 0} (\lambda_0\lambda_1)^l (\lambda_1\lambda_1)^r \right) \lambda_1\lambda_0 
\end{align*}
and the same reasoning gives for the probability of moving to the left
\begin{align*}
P(X(\tau_{n+1})=4m-4|\, X(\tau_{n})=4m) = \lambda_1\lambda_1 \left( \sum_{l,r\geq 0} (\lambda_0\lambda_1)^l (\lambda_1\lambda_1)^r \right) \lambda_0\lambda_0 ,
\end{align*}
which shows that the process $X(\tau_n)$ is a simple symmetric random walk with holding, and therefore visits the origin infinitely often.

For $N=3$, transience of $R$ is proven in Lemma 5 in \cite{Henry3colours}. If $N>3$, $\mathcal{T}$ contains a subtree on which $R$ is transient, so $R$ is transient on $\mathcal{T}$ as well. \qed

Finally, we give a statement and an open question for the case $d\geq 2$.
In order to make sure that the RWRE visits infinitely many sites, assume that $\mu$ is concentrated on the subset
$\widetilde{\cal P} ^d = \{\gamma \in {\cal P}^d: \gamma(e_i) > 0, \gamma(-e_i) > 0, 1 \leq i \leq d\}$, and let 
$\mathcal{M} ^d$ be the set of probability measures on $\widetilde{\cal P} ^d$.

\begin{thm}\label{thm2}
Assume $d\geq 2$ and assume that $\mu$ has a non-atomic part. Then, there exists a measurable mapping 
$\mathcal{A}:(\widetilde{\cal P} ^d)^{\mathbb{N}_0}\rightarrow \mathcal{M}^d$, such that for any measure $\mu \in \mathcal{M}^d$
\begin{align*}
P_\omega \big( \mathcal{A}(\xi) = \mu \big) =1 
\end{align*}
for $P$-almost all $\omega$.
\end{thm}
The proof of Theorem \ref{thm2} goes along the same lines as the informal description after Theorem \ref{haupt1}, which showed how to reconstruct $\mu$ in the case where $\mu$ has a non-atomic part.

\begin{quest}\label{quest}
Assume $d\geq 2$ and assume that $\mu$ is purely atomic. Is there a measurable mapping 
$\mathcal{A}:(\widetilde{\cal P} ^d)^{\mathbb{N}_0}\rightarrow \mathcal{M}^d$, such that for any measure $\mu \in \mathcal{M}^d$
\begin{align*}
P_\omega \big( \mathcal{A}(\xi) = \mu \big) =1
\end{align*}
for $P$-almost all $\omega$?
\end{quest}

\textbf{Acknowledgement:} We thank Noam Berger for discussions.

\bibliographystyle{alpha}

\begin{thebibliography}{[PE08]}

\bibitem{AdelmanEn} Adelman, O. and Enriquez, N. (2004)
\textit{Random walks in random environment: what a single trajectory tells.}
Israel. J. Math. \textbf{142}, 205-220.

\bibitem{BenKe} Benjamini, I. and Kesten, H. (1996)
\textit{Distinguishing sceneries by observing the scenery along a random walk path.}
J. Anal. Math. \textbf{69}, 97-135.

\bibitem{FrankJeffsurvey} Den Hollander, F. and Steif, J. (2006)
\textit{Random walk in random scenery: a survey of some recent results.}
Dynamics and stochastics, IMS Lecture Notes Monogr. Ser., \textbf{48}, Inst. Math. Statist., Beachwood, OH, 53-65.
							
\bibitem{Kestensurvey} Kesten, H. (1998)
\textit{Distinguishing and reconstructing sceneries from observations along random walk paths.}
Microsurveys in discrete probability, DIMACS Ser. Discrete Math. Theoret. Comput. Sci., \textbf{41}, Amer. Math. Soc., Providence, RI, 75-83.

\bibitem{Henry3colours} Matzinger, H. (1999)
\textit{Reconstructing a three-color scenery by observing it along a simple random walk path.}
Random Structures Algorithms \textbf{15}, no 2, 196-207.

\bibitem{Henry2colours} Matzinger, H. (2005)
\textit{Reconstructing a two-color scenery by observing it along a simple random walk path.}
Ann. Appl. Probab. \textbf{15}, 778–-819. 

\bibitem{HenryMatthias} L\"owe, M. and Matzinger, H. (2002).
\textit{Scenery reconstruction in two dimensions with many colors.}
Ann. Appl. Probab. \textbf{12} no. 4, 1322–-1347.

\bibitem{So} Solomon, F.(1975) 
\textit{Random walks in random environments.}
Ann. Probab. \textbf{3}, 1-31. 		

\bibitem{Z}  Zeitouni, O. (2004) \textit{Random walks in random environment.}, 
Lecture Lectures on probability theory and
		 	statistics, Lecture Notes in Math. \textbf{1837}, Springer, Berlin, 189-312.	


\end{thebibliography}

\bigskip
{\footnotesize

Nina Gantert: Technische Universit\"at M\"unchen,
Fakult\"at f\"ur Mathematik,
Boltzmannstra\ss e~3, 
85748~Garching bei M\"unchen,
Germany, gantert@ma.tum.de\\

Jan Nagel: Technische Universit\"at M\"unchen,
Fakult\"at f\"ur Mathematik,
Boltzmannstra\ss e~3, 
85748~Garching bei M\"unchen,
Germany, jan.nagel@ma.tum.de\\
}

\end{document}